\documentclass[graybox]{svmult}

\usepackage{type1cm}        
\usepackage{makeidx}         
\usepackage{graphicx}        
\usepackage{multicol}        
\usepackage[bottom]{footmisc}

\usepackage{newtxtext}       %
\usepackage[varvw]{newtxmath}       

\usepackage[colorlinks=true, allcolors=blue]{hyperref}

\makeindex   

\begin{document}
	
\title*{Schwarz Modulus Based Matrix Splittings with Minimal Polynomial Extrapolation Acceleration for linear complementarity problems arising from American option pricing} 
\titlerunning{The MPE method for the LCP arising from American option pricing}

\author{Martin J. Gander\orcidID{0000-0001-8450-9223} and Si-Wei Liao\orcidID{0000-0001-5592-6597} and Liu-Di Lu\orcidID{0000-0003-3629-9879}}
\institute{Martin J. Gander \at University of Geneva, Switzerland, \email{martin.gander@unige.ch}
\and Si-Wei Liao \at University of Geneva, Switzerland and Lanzhou University, China, \email{siwei.liao@unige.ch}
\and Liu-Di Lu \at Lund University, Sweden, \email{liudi.lu@math.lu.se}}
	%
	%
\maketitle
	
\abstract*{Pricing American options is more complicated than
          pricing European options, because they can be exercised at
          any time, and one thus needs to solve a linear
          complementarity problem instead of simply doing time
          stepping for computing European options. We introduce a new
          Schwarz modulus-based splitting method for solving such
          linear complementarity problems, and further accelerate them
          using Modified Polynomial Extrapolation, a non-linear vector
          sequence acceleration technique, which is very much related
          to Krylov methods in the linear case. Numerical experiments
          on a model problem show that our new solver can have close
          to an order of magnitude lower iteration counts than the
          classically used modulus-based matrix splitting technique.}

	\section{Introduction}

	Options are one of the most fundamental financial derivative
	instruments in modern markets. They give the holder the right,
	but not the obligation, to buy (call option) or sell (put option) an
	underlying asset at a fixed strike price, either on or before a
	specified expiration date. The evaluation of options has long
	been an active research field. One of the most famous models for
	option pricing is the Black--Scholes model, introduced in 1973 by
	Black and Scholes~\cite{BlackandScholes1973The}. Based on when the
		options can be exercised, they are classified as European and
	American options. European options can only be exercised at
	expiration, and their prices can be obtained analytically by solving
	the Black--Scholes equation, which means solving a diffusive
		boundary value problem. In contrast, American options may be
	exercised at any time before expiration. This early exercise feature
	imposes an additional constraint on the option value, transforming the
	Black--Scholes model into a free boundary value problem. Due to this
	complexity, American options generally do not have closed-form
	solutions, and numerical methods must be used.
	
	A simple numerical method for pricing American options is the
	binomial method proposed by Cox, Ross and
	Rubinstein~\cite{CoxRoss1979Option}. It converges however slowly
	and is difficult to extend to higher dimensions. In contrast, the
	finite difference method
	(FDM)~\cite{BrennanSchwartz1977The,BrennanSchwartz1978Finite} is
	widely used, offering both high accuracy and extensibility to
	higher dimensions. Discretizing the Black--Scholes model in both
	space and time, the FDM transforms the continuous partial
	differential equation (PDE) into a discrete algebraic system. For
	American options, the resulting discretized model can be reformulated
	as a linear complementarity problem (LCP)~\cite{HanWu2004A}, whose
	efficient numerical solution is an active research topic.
	Recently, the modulus-based matrix splitting (MMS)
        method~\cite{Bai2010Modulus} has been introduced as an
        efficient approach for solving LCPs, due to good numerical
        properties and convergence. One first uses a transformation to
        reformulate the LCP into an equivalent implicit fixed-point
        equation, and then applies a modulus-based matrix splitting
        strategy to construct an iterative scheme. Compared with
        conventional methods, such as the projected successive over
        relaxation method~\cite{Cryer1971PSOR} and the fixed-point
        method in~\cite{ShiYang2016Fixed}, modulus based methods
        avoid projection operations at each iteration, reducing
        computational cost.
	
	To improve the convergence of MMS for American options
        pricing, we introduce in this paper a new Schwarz modulus
          based matrix splitting method, and the minimal polynomial
        extrapolation (MPE), a non-linear vector sequence
          convergence acceleration technique, which is very much
          related to Krylov methods in the linear case
          \cite{McCoid2023Extrapolation}. In Section~\ref{Model-AOP},
        we reformulate the Black--Scholes model as a standard LCP
        using variable transforms and a finite difference
        discretization. We present in
        Section~\ref{numerical-method} the new Schwarz modulus
          based splitting method and its MPE acceleration, followed in
          Section \ref{numerical-test} by numerical results, and some
        concluding remarks.
	
	\section{Model, transformation and discretization of American options} \label{Model-AOP}
	
	We denote by $V(S,t)$ the price of an American option and by
	$G(S)$ the given payoff function. Here, $S>0$ is the underlying asset
	value and $t\in[0,T]$ is the time variable with $T$ the expiration
	time. The Black--Scholes operator $\mathcal{L}$ is defined as 
	$\mathcal{L} := \frac{\partial}{\partial t} + \frac{1}{2} \sigma^2 S^2 \frac{\partial^2}{\partial S^2} + (r-\delta)S \frac{\partial}{\partial S} - r$,
	where $\sigma$ is the fixed volatility, $r$ is the constant risk free
	interest rate, and $\delta$ is the continuous dividend yield. Due
	to the early exercise possibility of American options, $V(S,t)$ should
	be larger than $G(S)$ in order to avoid arbitrage possibilities. Based
	on a standard no arbitrage argument, $V(S,t)$ satisfies the 
	complementarity problem
	\begin{equation}\label{BS-model}
		\mathcal{L}V\leq0, \quad 
		V(S,t)\geq G(S), 
		\quad \text{and} \quad 
		\mathcal{L}V\cdot\left(V-G\right)=0,
	\end{equation}
	with the final condition $V(S,T) = G(S)$. For a fixed strike price
	$K$, the payoff function $G(S)$ is defined as
\begin{equation}\label{payoff-function}
\text{call option: }\ G(S) =\max\{S - K, 0\},
\quad
\text{put option: }\ G(S) =\max\{K-S, 0\}.
\end{equation}
	To complete the problem, we have the boundary conditions
	$\lim\limits_{S\rightarrow0}V(S,t) = K$ and
	$\lim\limits_{S\rightarrow\infty}V(S,t) = 0$ for put options, and
	$\lim\limits_{S\rightarrow0}V(S,t) = 0$ and
	$\lim\limits_{S\rightarrow\infty}V(S,t) \sim S - K\exp(-r (T-t))$ for
	call options.
	Note that the asset value $S$ in practice is in a range
	of $(S_{\min}, S_{\max})$, where typically Dirichlet boundary
	conditions are imposed.
	
	A change of variables can transform the Black--Scholes
	equation into a standard heat equation. This transformation simplifies
	the structure of the Black--Scholes equation, and improves the stability properties for numerical methods. As shown
	in~\cite{TavellaRandall2000Pricing}, this transformation leads to a
	flattened eigenvalue distribution of the discretized PDE
	operator. Hence, numerical schemes applied to the transformed
	heat equation achieve better stability behavior compared to those
	applied directly to the original Black--Scholes equation. For this
	reason, we introduce the change of variables
	\begin{equation}\label{Variable-Transformation}
		t = T - \dfrac{2\tau}{\sigma^2}, \quad 
		S = K e^x, \quad \text{and}\quad 
		V(S,t) = K e^{\alpha x + \beta \tau} u(x,\tau),
	\end{equation}
	where $\alpha = -(h_\delta - 1)/2$, $\beta = -(h_\delta - 1)^2/4 - h$,
	$h = 2r/\sigma^2$, and $h_\delta = 2(r - \delta)/\sigma^2$. The
	original time interval $t\in [0,T]$ now becomes
	$\tau\in[0,\sigma^2T/2]$, where the new variable $\tau$
	propagates backward with respect to $t$, and represents the remaining
	lifetime of the American option. For the asset value $S$, the
	logarithmic transformation maps the original semi-infinite domain
	$S\in (0,\infty)$ to the entire real line $x\in
	(-\infty,\infty)$. Using the change of
	variables~\eqref{Variable-Transformation}, we can then derive a
	complementarity problem for the new unknown $u(x,\tau)$,
	\begin{equation}\label{BSmodel-heatEquation}
		-\frac{\partial u}{\partial \tau} + \frac{\partial^2 u}{\partial x^2} \leq 0, 
		\quad u(x,\tau) \geq g(x,\tau), \quad \text{and}\quad 
		\left( -\frac{\partial u}{\partial \tau} + \frac{\partial^2 u}{\partial x^2} \right)(u - g) = 0,
	\end{equation}
	which is equivalent to the Black--Scholes complementarity
	problem~\eqref{BS-model}. Furthermore, the payoff
	function~\eqref{payoff-function} now becomes
	\begin{equation}\label{transformed-payoff-function}
		g(x,\tau)=\begin{cases}
			e^{\frac{\tau}{4}\left((h_\delta - 1)^2 + 4h\right)} \max\left\{ e^{\frac{x}{2}(h_\delta + 1)} - e^{\frac{x}{2}(h_\delta - 1)},\ 0 \right\}, & \text{for a call option}, \\
			e^{\frac{\tau}{4}\left((h_\delta - 1)^2 + 4h\right)} \max\left\{ e^{\frac{x}{2}(h_\delta - 1)} - e^{\frac{x}{2}(h_\delta + 1)},\ 0 \right\}, & \text{for a put option}.
		\end{cases}
	\end{equation}
	The final condition $V(S,T) = G(S)$ becomes an initial condition
	$u(x,0) = g(x,0)$, and the boundary conditions of $V(S,t)$ are
	now $\lim\limits_{x\rightarrow \pm\infty} u(x,\tau) =
	\lim\limits_{x\rightarrow \pm\infty} g(x,\tau)$. Once we have solved for
	$u(x,\tau)$, we can recover the original solution $V(S,t)$ using
	the change of variables~\eqref{Variable-Transformation}.
	
	To apply the FDM to discretize the transformed complementarity
        problem~\eqref{BSmodel-heatEquation}, we replace the infinite
        spatial domain $(-\infty, \infty)$ by a truncated interval $(a, b)$ with $a<0<b$. 
		Such a truncated interval is chosen based on the economical relevance of the asset value $S$, e.g., prices beyond which the option value is negligible. We denote this range of the asset value by $(S_{\min},S_{\max})$, and choose $a, b$ such that $Ke^a \leq S_{\min}, Ke^b \geq S_{\max}$ to guarantee numerical accuracy within $(S_{\min},S_{\max})$.
        Note that the choice of the interval has
        a direct impact on the accuracy of the numerical
        approximation. In particular, $a$ should be chosen
        sufficiently small and $b$ sufficiently large to keep the
        truncation error within an acceptable tolerance. However,
        taking $a$ too small or $b$ too large increases computational
        cost. Thus, there is a trade-off between the numerical
        accuracy and the computational cost. One could also use
          absorbing boundary conditions for better truncation
          properties, but this is not well explored; practical
        choices of the truncation interval can be found, 
        e.g., in~\cite{KangroNicolaides2000Far}.

	We apply a uniform grid to discretize the space-time domain
	$(a,b)\times (0,\sigma^2T/2]$. Let $\Delta x = \frac{b -
		a}{n}$ and $\Delta \tau = \sigma^2 T/(2m)$ denote the mesh
	sizes in $x$ and $\tau$ direction. Define the grid points as
	$x_i=a+i\Delta x$, $i=0, \dots, n$, and $\tau_j=j\Delta \tau$, $j=0,
	\dots, m$, and denote the approximate solution by $u_i^j \approx
	u(x_i, \tau_j)$. A $\theta$-weighted FDM applied to the first
	inequality in~\eqref{BSmodel-heatEquation} gives
$u_i^{j+1} - \lambda \theta (u_{i+1}^{j+1} - 2u_i^{j+1} + u_{i-1}^{j+1}) \geq u_i^j + \lambda (1 - \theta)(u_{i+1}^j - 2u_i^j + u_{i-1}^j)$,
	for $i=1, \dots, n-1$, $j=0, \dots, m-1$, $\lambda := \Delta \tau / \Delta x^2$, and $\theta\in(0, 1)$.
          In the matrix-vector form, we have
        $A\boldsymbol{u}^{j+1} \geq B \boldsymbol{u}^{j} +
        \boldsymbol{b}^{j+1} =: \boldsymbol{r}^j$, where $A = I_{n-1}
        + \lambda \theta T$, $B = I_{n-1} - \lambda (1-\theta) T$ and
        $\boldsymbol{u}^{j}=[u_1^j, \dots, u_{n-1}^j]^\top$. Here,
        $I_{n-1}\in \mathbb{R}^{(n-1)\times(n-1)}$ denotes the
        identity matrix, $T\in \mathbb{R}^{(n-1)\times(n-1)}$ is a
        tridiagonal matrix and $\boldsymbol{b}^{j+1}\in \mathbb{R}^{n-1}$
        is a vector, namely
	\[
	T=\begin{bmatrix}
		2 & -1 &        &        &        \\
		-1 & 2 & -1     &        &        \\
		& \ddots & \ddots & \ddots &   \\
		&        & -1     & 2 & -1 \\
		&        &        & -1 & 2
	\end{bmatrix},\qquad
	\boldsymbol{b}^{j+1} =
	\begin{bmatrix}
		\lambda(1{-}\theta) g_0^j + \lambda\theta g_0^{j+1} \\
		0 \\ 
		\vdots \\ 
		0 \\
		\lambda(1{-}\theta) g_n^j + \lambda\theta g_n^{j+1}
	\end{bmatrix},
	\]
	with $g_i^j \approx g(x_i, \tau_j)$. Then, the FDM discretized
	complementarity problem~\eqref{BSmodel-heatEquation} becomes a
	sequence of complementarity problems:
$A\boldsymbol{u}^{j+1} - \boldsymbol{r}^j \geq 0$, $\boldsymbol{u}^{j+1} \geq \boldsymbol{g}^{j+1}$, and $(A \boldsymbol{u}^{j+1} - \boldsymbol{r}^{j})^\top (\boldsymbol{u}^{j+1} - \boldsymbol{g}^{j+1}) = 0$,
	with the initial condition $u_i^0 = g_i^0$, $i=1, \cdots, n-1$, and
	boundary conditions $u_0^j = g_0^j$, $u_n^j = g_n^j$, $j=0, \cdots,
	m-1$. For computational convenience, we can
	write it in the standard form
	of an LCP. To achieve this, we define $\boldsymbol{z} :=
	\boldsymbol{u}^{j+1} - \boldsymbol{g}^{j+1}$ and $\boldsymbol{q} := A
	\boldsymbol{g}^{j+1} - \boldsymbol{r}^j$ and re-write the early
	exercise constraint $\boldsymbol{u}^{j+1} \geq \boldsymbol{g}^{j+1}$
	as a nonnegativity condition $\boldsymbol{u}^{j+1} -
	\boldsymbol{g}^{j+1} \geq 0$. This gives the standard LCP form
	\begin{equation}\label{LCP}
		A \boldsymbol{z} + \boldsymbol{q} \geq 0, \quad \boldsymbol{z} \geq 0 
		\quad \text{and} \quad 
		(A \boldsymbol{z} + \boldsymbol{q})^\top \boldsymbol{z} = 0.
	\end{equation}
	Note that the vector $\boldsymbol{q}$ is given, since
	$\boldsymbol{g}^{j+1}$ is given and $\boldsymbol{r}^j$ is known at the
	$j$th time step. The only unknown now is $\boldsymbol{z}$, which is
	the solution of the LCP~\eqref{LCP}. In other words, one needs to
	solve an LCP at each time step, which motivates the development of
	efficient numerical methods.
	
\section{Schwarz modulus based matrix splitting and minimal polynomial extrapolation}\label{numerical-method}
	
We now present numerical methods for solving the LCP~\eqref{LCP}. The
classical MMS consists of applying a modulus transformation to
preserve both the complementarity and nonnegativity conditions
in~\eqref{LCP}. It has been shown in~\cite[Theorem
  2.1]{Bai2010Modulus} that the LCP~\eqref{LCP} is equivalent to the
fixed point problem: find $\mathbf{y}\in\mathbb{R}^{n-1}$ such that
\begin{equation}\label{fixed-point-equation}
  (M+\Omega)\mathbf{y}=N\mathbf{y}+(\Omega-A)|\mathbf{y}|-\eta\boldsymbol{q},
\end{equation}
where $A=M-N$ is a matrix splitting of $A$, typically Gauss-Seidel,
$\Omega\in\mathbb{R}^{(n-1)\times (n-1)}$ is a positive diagonal
matrix, and $\eta>0$ is a constant. A detailed discussion on the choices of $\eta$ and $\Omega$ is given in~\cite[Section 3]{Bai2010Modulus}.
In particular, if $\mathbf{y}$
satisfies~\eqref{fixed-point-equation},
$\boldsymbol{z}=(|\mathbf{y}|+\mathbf{y})/\eta$ solves the
LCP~\eqref{LCP}. On the other hand, if $\boldsymbol{z}$ is a solution
of the LCP \eqref{LCP}, then
$\mathbf{y}=\eta\Omega^{-1}((\Omega-A)\boldsymbol{z}-\boldsymbol{q})/2$
satisfies the fixed-point problem~\eqref{fixed-point-equation}. Hence,
one can solve the fixed point problem~\eqref{fixed-point-equation} and
then recover the solution for the LCP~\eqref{LCP} using their
equivalence. A natural stationary iteration to solve~\eqref{fixed-point-equation} is given by $(M+\Omega)\mathbf{y}_{k+1}=N\mathbf{y}_{k}+(\Omega-A)|\mathbf{y}_{k}|-\eta \boldsymbol{q}$,
with a given initial vector $\mathbf{y}_{0}\in\mathbb{R}^{n-1}$, and
the iteration index $k=0,1,\ldots$. Once $\mathbf{y}_{k+1} $ is
computed, the corresponding approximate solution of the original
LCP~\eqref{LCP} can be recovered using
$\boldsymbol{z}_{k+1}=(|\mathbf{y}_{k+1}|+\mathbf{y}_{k+1})/\eta$.
		
We introduce now as a new solver an alternating Schwarz method: we
  decompose the spatial domain $(a, b)$ into subdomains with minimal
overlap of one mesh size $\Delta x$, which corresponds to a
block Gauss--Seidel method, see e.g.~\cite{gander2008schwarz}. This decomposition permits an arbitrary number of overlapping subdomains. For our problem, we observe from our numerical experiments that a two-subdomain decomposition gives the best performance.
So we focus here on two subdomains, which gives
\[
A=\begin{pmatrix} A_{11} & A_{12} \\ A_{21} & A_{22} \\ \end{pmatrix},\,\,\, M = \begin{pmatrix} A_{11} & 0 \\ A_{21} & A_{22} \\ \end{pmatrix}\,\,\,\text{and}\,\,\, N = \begin{pmatrix} 0 & -A_{12} \\ 0 & 0 \\ \end{pmatrix},
\]
where $A_{11}$ is the matrix corresponding to the discrete problem in
the spatial subdomain $(a, \Delta x)$ to the left of the origin
  (recall that $a<0$) with overlap $\Delta x$, and $A_{22}$ is the
corresponding matrix in the subdomain $(0, b)$ to the right of the
  origin. The matrix $A_{12}$ corresponds to the Dirichlet
transmission condition at the interface $x = \Delta x$, similarly
$A_{21}$ is the matrix for the Dirichlet transmission condition at $x
= 0$. For comparison, we also consider the classical point-wise
Gauss--Seidel splitting with $M=D+L$ and $N=-U$, where $D$, $L$ and
$U$ are the diagonal, strictly lower triangular, and strictly upper
triangular matrices.
	
Note that we need to solve at each time step an LCP
problem~\eqref{LCP}, and when there are many time steps, the
  total iteration number can be very large to solve all the
fixed point problems~\eqref{fixed-point-equation}. To further
accelerate this process, we consider Modified Polynomial
  Extrapolation (MPE) introduced by Cabay and
Jackson~\cite{CabayJackson1976Polynomial}, which is a powerful
acceleration technique for improving the convergence of
non-linear vector sequences. It is very effective for
accelerating fixed-point iterations used to solve linear or nonlinear
systems of equations, especially those arising in the discrete
solution of continuum problems, and is related to Krylov methods
  in the linear case, see e.g. \cite[Section 11.6 and 11.7]{gander2014scientific} and \cite{McCoid2023Extrapolation}.

For a sequence of vectors $\mathbf{y}_{0}, \mathbf{y}_{1},
\mathbf{y}_{2}, \cdots, \mathbf{y}_{k+1}$ generated by the MMS method,
we first define the sequence of the difference between two vectors as
$\boldsymbol{u}_i := \mathbf{y}_{i+1} - \mathbf{y}_i$ for $i = 0,
1, \ldots, k$, and put them into the $(n-1)\times(k+1)$ matrices
  $U_k := (\boldsymbol{u}_0, \boldsymbol{u}_1, \cdots,
  \boldsymbol{u}_k)$. The MPE approximation $\mathbf{s}_{k}$ is given
by the linear combination $\mathbf{s}_k = \sum_{i=0}^k \gamma_i
\mathbf{y}_i$, where the coefficients $\gamma_i$ are determined such
that $U_k \boldsymbol{\gamma} \approx 0$ and $\sum_{i=0}^{k}\gamma_i =
1$, computed in two steps: First, solve the linear system $U_{k-1} \boldsymbol{c} \approx
    -\boldsymbol{u}_k$ as a least squares problem, where
    $\boldsymbol{c} = (c_0, c_1, \dots, c_{k-1})^{T}$.
Second, set $c_k = 1$, and compute the scaled coefficients by
    $\gamma_i = c_i /\sum_{j=0}^{k} c_j$ for $i = 0, 1, \ldots,
    k$. When $\sum_{j=0}^{k} c_j = 0$, the extrapolated vector
    $\mathbf{s}_k$ does not exist.

\section{Numerical experiments and discussion}\label{numerical-test}

\begin{figure}[t]
  \centering
    \mbox{
  \includegraphics[width=0.45\textwidth]{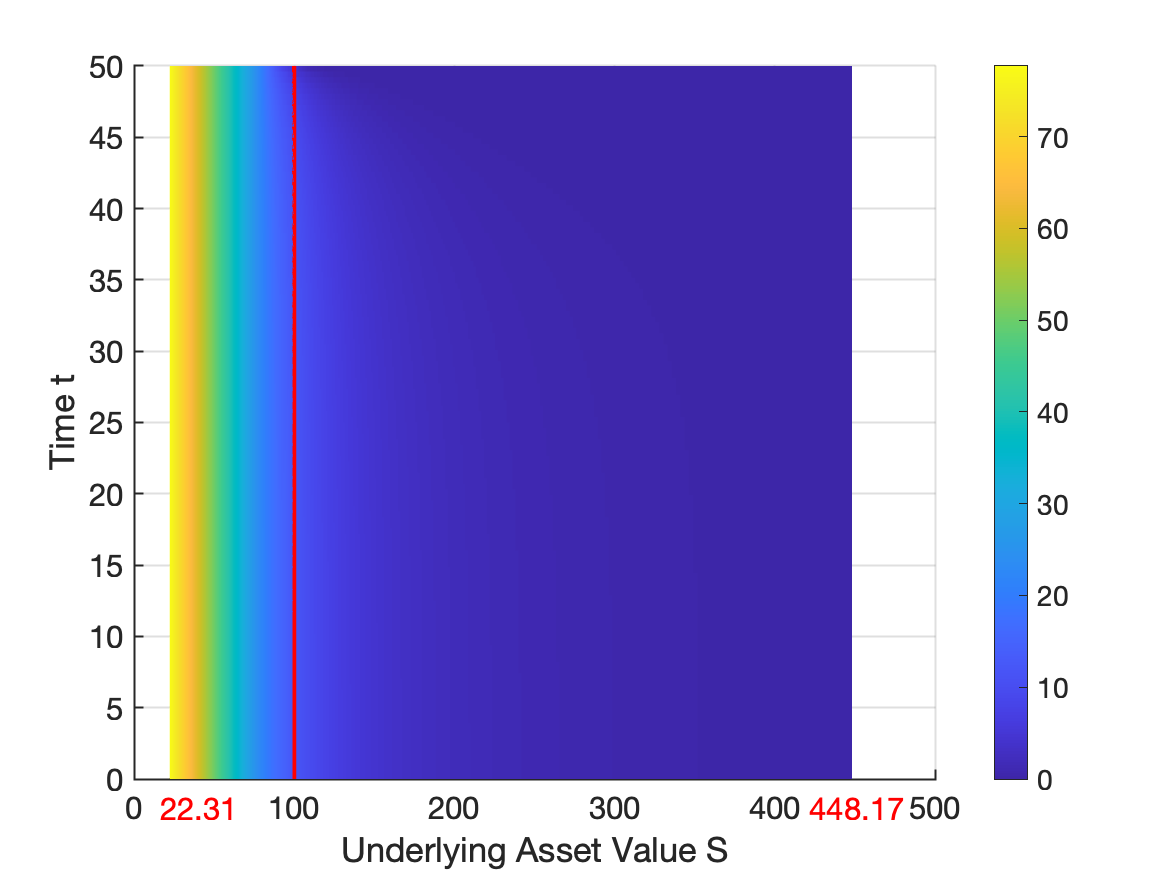}\quad
  \includegraphics[width=0.45\textwidth]{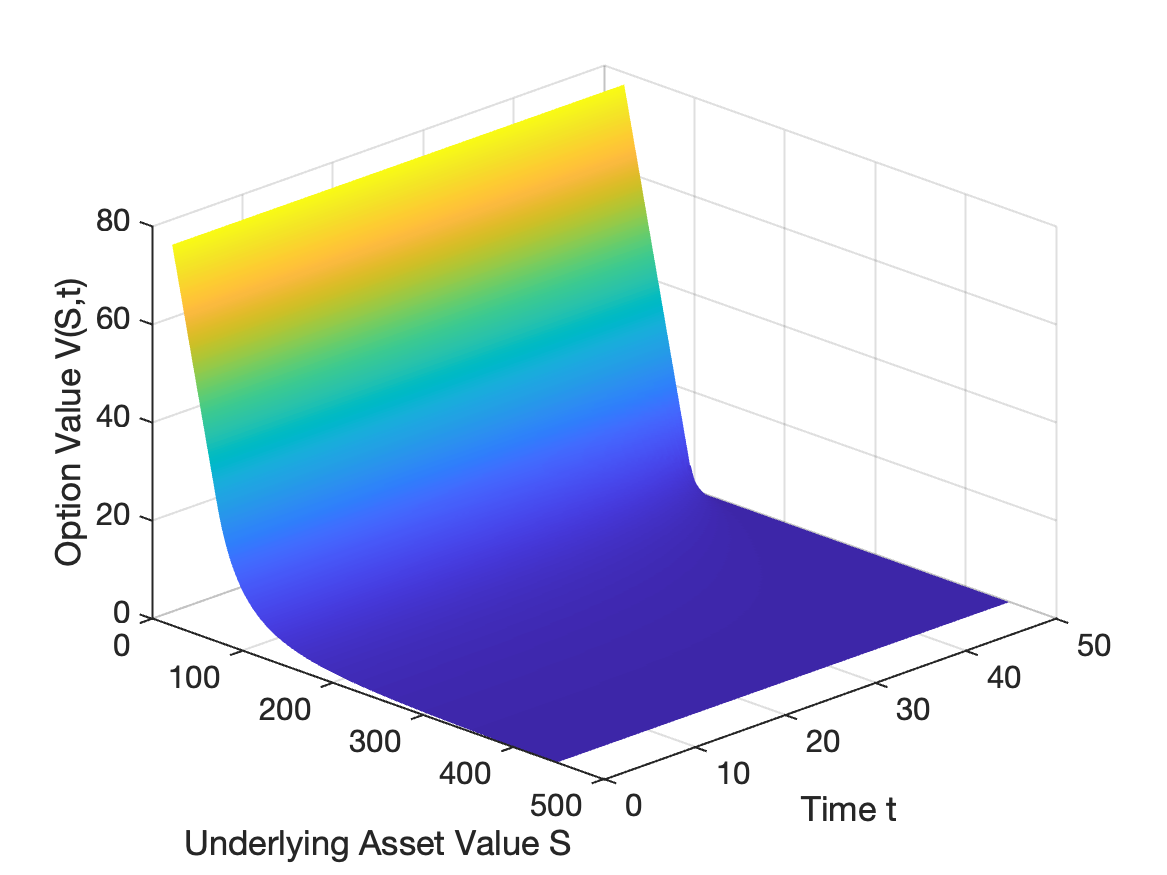}}
  \caption{Option value $V$ as a function of the asset value $S$ and
    time $t$ using our new Schwarz splitting with a mesh size
    $(\Delta x,\Delta \tau)=(2^{-7}, 2^{-7})$. Left: 2D view. Right:
    3D view.}
  \label{fig:option_surface}
\end{figure}

We illustrate now the numerical performance of our new
  Schwarz MMS and the new MPE acceleration for solving the
LCP arising from American option pricing. We consider an American put
option on a single underlying asset $S$. We set $K=100$, $r=0.05$,
$\sigma=0.2$, $\delta=0$, and $T=50$. For the FDM discretization, we
set $(a,b)=(-1.5,1.5)$ and use the Crank--Nicolson scheme with $\theta
=1/2$. 
We show in
Fig.~\ref{fig:option_surface}
the American put option price $V$ as a function of the underlying
asset value $S$ and time $t$, computed using our new Schwarz
MPE. The asset value ranges from $S_{\min}=Ke^a\approx22.31$ to
$S_{\max}= Ke^b\approx448.17$. Within this interval $(S_{\min},
  S_{\max})$, the option value $V$ can be determined at any time
$t$. In Fig.~\ref{fig:option_surface} (left), the option value $V$
decreases steadily as the asset value $S$ increases, consistent with
classical option pricing theory. Given the strike price $K=100$, when
$S$ is significantly greater than $K$, the option value $V$ approaches
zero. Fig.~\ref{fig:option_surface} (right) further shows that as
expiration approaches, the put option value curve becomes steeper with
increasing $S$, emphasizing the accelerated decay of the option value
when $S > K$.

In the MMS, we set $\eta=2$ and
$\Omega=1/2\text{diag}(A)$, and we test the classical point-wise Gauss--Seidel splitting, denoted by GS, two MPE accelerations, MPE-GS and MPECycle-GS, and the new Schwarz or Block Gauss-Seidel splitting,
denoted by BGS, with their MPE accelerated versions, MPE-BGS and MPECycle-BGS. In the first time step, the initial
vector $\mathbf{y}_{0}$ is set to the prescribed initial condition,
{\it i.e.}, all entries are initialized to $0$. For subsequent time
steps, the initial vector is taken as the numerical solution obtained
from the previous time step. At each time step, the stationary
iteration is terminated when the residual satisfies $\|
A(|\mathbf{y}_{k+1}|+\mathbf{y}_{k+1})+\Omega(\mathbf{y}_{k+1}-|\mathbf{y}_{k+1}|)+
\eta \boldsymbol{q} \| _{\infty}<10^{-6}$.

To compare the performance of the classical MMS and new BGS
  splitting with and without MPE acceleration, we show in
Table~\ref{tab:1}
the iteration numbers and computational time in seconds obtained with different mesh sizes for a
  convergence tolerance of $10^{-6}$. For each method, we show the total CPU time, total
number of iterations over all time steps and the average number
of iterations, which is computed by dividing the total iterations
over all time steps by the number of time steps.
\begin{table}
  \setlength{\tabcolsep}{1.2em}
  \centering
  \caption{Comparison of our new Schwarz splitting strategy (BGS)
    to classical Gauss-Seidel (GS) and their MPE
    accelerated versions with different mesh sizes}
  \label{tab:1}
  \begin{tabular}{l c c c c}
			\hline
			$(\Delta x, \Delta \tau)$ & Method &  Total Iterations & Average Iterations & CPU Time \\
			\hline
			($2^{-4}$, $2^{-4}$)   & GS & 1403 & 87  & 0.023839\\
			& MPE-GS & 429 & 27& 0.021326\\
			& MPECycle-GS & 480 & 30 & 0.007024 \\
			& BGS & 865 & 54 & 0.016162\\
			& MPE-BGS & 297 & 19 & 0.012260 \\
			& MPECycle-BGS & 363 & 23 & 0.005096 \\
			\hline
			($2^{-5}$, $2^{-5}$)   & GS & 5491 & 171 & 0.072522\\
			& MPE-GS & 1753& 55 & 0.264530 \\
			& MPECycle-GS & 1746 & 55 & 0.027914 \\
			& BGS & 3168 & 99 & 0.058267 \\
			& MPE-BGS & 936 & 29 & 0.066916 \\
			& MPECycle-BGS & 1072 & 34 & 0.022632\\
			\hline
			($2^{-6}$, $2^{-6}$)   & GS & 21123 & 330 & 0.619466\\
			& MPE-GS & 6542 & 102 & 5.260333\\
			&MPECycle-GS & 6368 & 100 & 0.230351\\
			& BGS & 11666 & 182 &0.418613 \\
			& MPE-BGS & 2887 & 45 & 0.443999\\
			&MPECycle-BGS & 3328 & 52& 0.144547\\
			\hline
		\end{tabular}
	\end{table}	
For the pointwise splitting case, we observe that compared with GS, MPE-GS divides the required iteration number for convergence by a factor of 3.21 on average. However, the performance in terms of the required CPU time is much different. As the mesh is refined, the required CPU time of MPE-GS increases drastically, due to the solving of the least square problem $U_{k-1} \boldsymbol{c} \approx -\boldsymbol{u}_k$ at each MPE iteration using the backslash of MATLAB\footnote{There are faster ways to implement the numerical solving of such least square problems by updating, like for Krylov methods, see e.g.,~\cite{Mouhssine2025}, but this is beyond the scope of our manuscript.}, and the proportion of CPU time needed rises sharply from 53.5\% to 92.2\%. Based on this observation, we apply an improved MPE acceleration method, MPECycle-GS, which consists of using the MPE acceleration only after every $N_c$ iterations of GS. This significantly improves the performance of the MPE acceleration, as it reduces both the frequency and size of solving least square problems, i.e., from solving for all previous iterates to only for the previous $N_c$ iterates. In our tests, $N_c=15$ gives a good improvement compared with GS, i.e., a speedup of 3.13 in terms of iteration numbers and 2.89 in terms of CPU time. This clearly shows that MPE is capable of a significant acceleration of the non-linear point-wise GS process.

We also see that our new Schwarz splitting variant BGS requires fewer iterations compared with GS. This is expected, since the Schwarz-based BGS takes advantage of the subdomain solves, yielding a stronger global correction at each iteration. The performance of both methods is similar in terms of CPU time. However, by
applying MPE to accelerate BGS, we see that the required CPU time of MPE-BGS is comparable to that of BGS, since the proportion of CPU time required for solving the least square problem is around 50\% on average, unlike a drastic increase in the point-wise case. When applying MPE every $N_c$ iterations, MPECycle-BGS can further improve the performance with a speedup of 2.95 in terms of iteration numbers and 2.88 in terms of CPU time compared with BGS. Overall, our novel Schwarz-based algorithms outperform the classical splitting methods, which motivates a further investigation of this new algorithm applied to LCP problems.

Due to the early exercise constraint, the American option pricing
problem is harder to solve numerically than the more classical
  European option pricing problem, and one must solve a sequence of
  linear complementarity problems (LCPs). To do so, we introduced a
  new Schwarz modulus bases splitting method for solving such LCPs,
  and further accelerated convergence using Modified Polynomial
  Extrapolation, a non-linear vector sequence acceleration
  technique related to Krylov methods in the linear case
  \cite{McCoid2023Extrapolation}.  Our new Schwarz modulus based provides a brand new avenue for numerically solving
the American option pricing problem.
	
\bibliographystyle{spmpsci}
\bibliography{sample}

@article{Bai2010Modulus,
title={Modulus--based matrix splitting iteration methods for linear complementarity problems},
author={Z. -Z. Bai},
pages={917--933},
year={2010},
volume={17},
journal={Numerical Linear Algebra with Applications}}

@article{BlackandScholes1973The,
title={The pricing of options and corporate liabilities},
author={F. Black and M. Scholes},
pages={637--659},
year={1973},
volume={81(3)},
journal={Journal of Political Economy}}

@article{BrennanSchwartz1977The,
title={The valuation of {A}merican put options},
author={M. J. Brennan and E. S. Schwartz},
pages={449--462},
year={1977},
volume={32(2)},
journal={Journal of Finance}}

@article{BrennanSchwartz1978Finite,
title={Finite Difference Methods and Jump Processes Arising in the Pricing of Contingent Claims: A Synthesis},
author={M. J. Brennan and E. S. Schwartz},
pages={461--474},
year={1978},
volume={13(3)},
journal={Journal of Financial and Quantitative Analysis}}

@article{CabayJackson1976Polynomial,
title={A polynomial extrapolation method for finding limits and antilimits of vector sequences},
author={S. Cabay and L. W. Jackson},
pages={734--752},
year={1976},
volume={13},
journal={SIAM Journal on Numerical Analysis}}

@article{CoxRoss1979Option,
title={Option pricing: A simplified approach},
author={J. C. Cox, S. A. Ross and M. Rubinstein},
pages={229--263},
year={1979},
volume={7(3)},
journal={Journal of Financial Economics}}

@article{Cryer1971PSOR,
title={The solution of a quadratic programming problem using systematic overrelaxation},
author={C. W. Cryer},
pages={385--392},
year={1971},
volume={9(3)},
journal={SIAM Journal on Control}}

@article{gander2008schwarz,
	title={Schwarz methods over the course of time},
	author={Gander, Martin J.},
	journal={Electron. Trans. Numer. Anal},
	volume={31},
	number={5},
	pages={228--255},
	year={2008}
}

@book{gander2014scientific,
  title={Scientific computing-An introduction using Maple and MATLAB},
  author={Gander, W. and Gander, M. J. and Kwok, F.},
  volume={11},
  year={2014},
  publisher={Springer Science \& Business}
}

@article{HanWu2004A,
title={A fast numerical method for the {B}lack--{S}choles equation of {A}merican options},
author={H. Han and X. -N. Wu},
pages={2081--2095},
year={2004},
volume={41(6)},
journal={SIAM Journal on Numerical Analysis}}

@article{KangroNicolaides2000Far,
title={Far field boundary conditions for {B}lack--{S}choles equations},
author={R. Kangro and R. Nicolaides},
pages={1357--1368},
year={2000},
volume={38},
journal={SIAM Journal on Numerical Analysis}}

@article{McCoid2023Extrapolation,
	title={Extrapolation methods as nonlinear {K}rylov subspace methods},
	author={McCoid, C. and Gander, Martin J.},
	journal={Linear Algebra and Its Applications},
	year={2023}}

@phdthesis{Mouhssine2025,
title = {Vector extrapolation methods with applications to geometric multigrid and nonlinear least-squares problems},
year = {2025},
author={Mouhssine, A.},
month={December},
school={Mohammed VI Polytechnic University}}

@article{ShiYang2016Fixed,
title={A fixed point method for the linear complementarity problem arising from {A}merican option pricing},
author={X. -J. Shi, L. Yang and Z. -H. Huang},
pages={921--932},
year={2016},
volume={32},
journal={Acta Math. Appl. Sin. Engl. Ser.}}

@book{TavellaRandall2000Pricing,
title={Pricing Financial Instruments: The finite difference method},
author={D. Tavella and C. Randall},
year={2000},
publisher={John Wiley and Sons, New York}}
	
\end{document}